\documentclass[titlepage,twoside,12pt]{article}
\usepackage{amssymb}
\usepackage{amsfonts}
\begin{document}

{\bf \Large Synthesizing Sums} \\ \\

Elem\'{e}r E Rosinger \\ \\
Department of Mathematics \\
and Applied Mathematics \\
University of Pretoria \\
Pretoria \\
0002 South Africa \\
eerosinger@hotmail.com \\ \\ \\

{\bf Abstract} \\

Polynomial functions $f : \mathbb{N}_+ \longrightarrow \mathbb{N}_+$ are studied for which
sums of arbitrary length $f ( 1 ) + f ( 2 ) + f ( 3 ) + \ldots + f ( n )$, with $n \in
\mathbb{N}_+$, can be expressed by polynomial functions $g : \mathbb{N}_+ \longrightarrow
\mathbb{N}_+$ which involve a bounded number of operations, thus not depending on $n$. Open
problems and extensions are presented. \\ \\ \\

{\bf 1. Remarks on Summation Formulas} \\

As is well known, we have the relations \\

(1.1)~~~ $ 1 + 2 + 3 + \ldots + n ~=~ n ( n + 1 ) / 2,~~~ n \in \mathbb{N}_+ $ \\

where $\mathbb{N}_+ = \{ 1, 2, 3, \ldots \}$. Similarly, we have \\

(1.2)~~~ $ 1^2 + 2^2 + 3^2 + \ldots + n^2 ~=~
                  n ( n + 1 ) ( 2 n + 1 ) / 6,~~~ n \in \mathbb{N}_+ $ \\

and \\

(1.3)~~~ $ 1 . 1 ! + 2 . 2 ! + 3 . 3 ! + \ldots + n . n ! ~=~
                                  ( n + 1 ) ! - 1,~~~ n \in \mathbb{N}_+ $ \\

On the other hand, as far as we know, no summation formula of the above type is available for
the simpler looking \\

(1.4)~~~ $ 1 ! + 2 ! + 3 ! + \ldots + n ! ~=~ ?,~~~ n \in \mathbb{N}_+ $ \\

{\bf Remarks} \\

1) It should be noted that the right hand terms in (1.1), (1.2) are polynomials in $n$, while
the left hand terms, as they stand, are not polynomials, since the number of terms they have
is dependent on $n$, and thus it is {\it not} bounded. \\
This difference between the left and right hand terms in (1.1), (1.2) is further highlighted
by the fact that the left hand terms in (1.1), (1.2), when considered as expressions of the
form \\

(1.4)~~~ $ f ( 1 ) + f ( 2 ) + f ( 3 ) + \ldots + f ( n ) $ \\

with $f : \mathbb{N}_+ \longrightarrow \mathbb{N}_+$, have the coefficients in $\mathbb{Z}$,
namely, these coefficients are each $1$, while the coefficients in the right hand terms in
(1.1), (1.2) are {\it no} longer in $\mathbb{Z}$, but in $\mathbb{Q}$, which is an {\it
extension} of $\mathbb{Z}$. \\

2) Similarly, the right hand term in (1.3) is a polynomial in $n$ and $n !$, while as they
stand, the left hand terms in (1.3), (1.4) are not such polynomials. \\

3) The nontriviality of relations such as (1.1), (1.2) can easily be seen in (A.1) in Lemma A
in the Appendix, which among other implies that, for every given $k \in \mathbb{N}_+$, the
{\it countably infinite} system of linear nonhomogeneous equation in $c_{k,\,0},~ c_{k,\, 1},~
\ldots ,~ c_{k,\, k} \in \mathbb{Q}$, namely \\

(1.5)~~~ $ n\, c_{k,\,0} + n^2\, c_{k,\, 1} + \ldots + n^{k + 1}\, c_{k,\, k} ~=~
                 1^k + 2^k + 3^k + \ldots + n^k,~~~ n \in \mathbb{N}_+ $ \\

has a unique solution. \\

4) The remark at 1) above about relations such as (1.1), (1.2), or more generally, those in
(2.6), (3.4) or (A.1) in the sequel, can also be seen {\it the other way round}. Namely, given
a ring $S$ and a subring $R$ of $S$, as well as an element $u \in R$, let us denote by \\

(1.6)~~~ $ S_{R,\,u}\, [ x ] $ \\

the set of all polynomials $g \in S\, [ x ]$, such that \\

(1.7)~~~ $ \begin{array}{l}
               \exists~~~ f \in R\, [ x ] ~:~ \\ \\
               \forall~~~ n \in \mathbb{N}_+ ~:~ \\ \\
               ~~~ g ( n . u ) ~=~ f ( u ) + f ( 2 . u ) + f ( 3 . u ) +
                                    \ldots + f ( n . u )
            \end{array} $ \\

\medskip
Then one can formulate the following \\

{\bf Open Problem 1.} \\

Give alternative characterizations for the polynomials $g$ in $S_{R,\,u}\, [ x ]$.

\hfill $\Box$ \\

Needless to say, one may similarly be interested in the relationship between the rings $R$ and
their ring extensions $S$, for which $S_{R,\,u}\, [ x ]$ is not void, for certain, or for all
$u \in R$. \\ \\

{\bf 2. Formulation of a Further Problem} \\

Let us denote by \\

(2.1)~~~ $ {\cal P}{\cal F}_\mathbb{Z} $ \\

the set of all functions $f : \mathbb{N}_+ \longrightarrow \mathbb{N}_+$, such that \\

(2.2)~~~ $ f ( n ) = p ( n, n ! ),~~~ n \in \mathbb{N}_+ $ \\

where $p ( x, y )$ is a two variable polynomial with coefficients in $\mathbb{Z}$. Further, we
denote by \\

(2.3)~~~ $ {\cal P}{\cal F}_\mathbb{Q} $ \\

the set of all functions $g : \mathbb{N}_+ \longrightarrow \mathbb{N}_+$, such that \\

(2.4)~~~ $ g ( n ) = q ( n, n ! ),~~~ n \in \mathbb{N}_+ $ \\

where $q ( x, y )$ is a two variable polynomial with coefficients in $\mathbb{Q}$. \\

Our interest is in the subset \\

(2.5)~~~ $ {\cal S}{\cal P}{\cal F}_\mathbb{Z} $ \\

of ${\cal P}{\cal F}_\mathbb{Z}$ which consists of all functions $f \in {\cal P}{\cal
F}_\mathbb{Z}$ which have the property \\

(2.6)~~~ $ \begin{array}{l}
                   \exists~~~ g \in {\cal P}{\cal F}_\mathbb{Q} ~:~ \\ \\
                   \forall~~~ n \in \mathbb{N}_+ ~:~ \\ \\
                   ~~~ f ( 1 ) + f ( 2 ) + f ( 3 ) + \ldots + f ( n ) ~=~ g ( n )
            \end{array} $ \\

Namely, we are looking for alternative {\it characterizations} of those functions $f \in {\cal
P}{\cal F}_\mathbb{Z}$ which belong to ${\cal S}{\cal P}{\cal F}_\mathbb{Z}$. \\

As for the non-triviality of this problem, we can note that, in view of (1.1) - (1.3), the
functions $f$ given respectively by \\

$~~~~~~ n \longmapsto n,~~~ n \longmapsto n^2,~~ n \longmapsto n . n ! $ \\

belong to ${\cal S}{\cal P}{\cal F}_\mathbb{Z}$, since they obviously satisfy (2.6), namely,
for the respective functions $g$ \\

$~~~~~~ n \longmapsto n ( n + 1 ) / 2,~~~
                    n ( n + 1 ) ( 2 n + 1 ) / 6,~~~ n \longmapsto ( n + 1 ) ! - 1 $ \\

On the other hand, as far as we know, the function $f$ given by \\

$~~~~~~ n \longmapsto n ! $ \\

may, or may not belong to ${\cal S}{\cal P}{\cal F}_\mathbb{Z}$. \\

Therefore, in view of the Theorem in the next section, we are led to \\

{\bf Open Problem 2.} \\

Does the equality \\

(2.7)~~~ $ {\cal S}{\cal P}{\cal F}_\mathbb{Z} ~=~ {\cal P}{\cal F}_\mathbb{Z} $ \\

hold ? \\ \\

{\bf 3. A Simpler Problem and its Solution} \\

We can limit ourselves in the above by not considering factorials, that is, by only
considering relations such as in (1.1), (1.2). In this case, let us denote by \\

(3.1)~~~ $ {\cal P}_\mathbb{Z} $ \\

the set of all functions $f : \mathbb{N}_+ \longrightarrow \mathbb{N}_+$, given by polynomials
with coefficients in $\mathbb{Z}$. Further, we denote by \\

(3.2)~~~ $ {\cal P}_\mathbb{Q} $ \\

the set of all functions $g : \mathbb{N}_+ \longrightarrow \mathbb{N}_+$, given by polynomials
with coefficients in $\mathbb{Q}$. Then, instead of ${\cal S}{\cal P}{\cal F}_\mathbb{Z}$ in
(2.5), we are considering its subset \\

(3.3)~~~ $ {\cal S}{\cal P}_\mathbb{Z} $ \\

of all functions $f \in {\cal P}_\mathbb{Z}$ which have the property \\

(3.4)~~~ $ \begin{array}{l}
                   \exists~~~ g \in {\cal P}_\mathbb{Q} ~:~ \\ \\
                   \forall~~~ n \in \mathbb{N}_+ ~:~ \\ \\
                   ~~~ f ( 1 ) + f ( 2 ) + f ( 3 ) + \ldots + f ( n ) ~=~ g ( n )
            \end{array} $ \\

and correspondingly, we are looking for alternative {\it characterizations} of those functions
$f \in {\cal P}_\mathbb{Z}$ which belong to ${\cal S}{\cal P}_\mathbb{Z}$. \\

Here however, in view of Lemma A in the Appendix, the situation is simple, since we have \\

{\bf Theorem} \\

(3.5)~~~ $ {\cal S}{\cal P}_\mathbb{Z} ~=~ {\cal P}_\mathbb{Z} $ \\

{\bf Proof.} \\

In view of (A.1), it follows that for every $k \in \mathbb{N}_+$, the function $n \longmapsto
n^k$ in ${\cal P}_\mathbb{Z}$ belongs to ${\cal S}{\cal P}_\mathbb{Z}$. \\

Let now be given any $f \in {\cal P}_\mathbb{Z}$, then \\

$~~~~~~ f ( n ) ~=~ a_0 + a_1\, n + a_2\, n^2 + \ldots + a_k\, n^k,~~~
                                                               n \in \mathbb{N}_+ $ \\

for a certain $k \in \mathbb{N}_+$ and suitable $a_0, a_1, a_2, \ldots , a_k \in
\mathbb{Z}$. \\

Let us now denote, see (A.1) in Lemma A in the Appendix \\

$~~~~~ \Sigma_n^k ~=~ 1^k + 2^k + 3^k + \ldots + n^k,~~~
         \Pi_n^k ~=~ c_{k,\,0}\, n + c_{k,\, 1}\, n^2 + \ldots + c_{k,\, k}\, n^{k + 1} $ \\

with $k, n \in \mathbb{N}_+$. Then \\

$~~~~~~ f ( 1 ) + f ( 2 ) + f ( 3 ) + \ldots + f ( n ) ~=~ $ \\

$~~~~~~ ~=~ a_0\, n + a_1\, \Sigma_n^1 + a_2\, \Sigma_n^2 +
                          \ldots + a_k\, \Sigma_n^k ~=~ $ \\

$~~~~~~ ~=~ a_0\, n + a_1\, \Pi_n^1 + a_2\, \Pi_n^2 +
                          \ldots + a_k\, \Pi_n^k,~~~ n \in \mathbb{N}_+ $ \\

and the last term is obviously a polynomial in $n$. \\ \\

{\bf 4. A General Scheme} \\

For every function $f : \mathbb{N}_+ \longrightarrow \mathbb{N}_+$, we define the function
$\Sigma f : \mathbb{N}_+ \longrightarrow \mathbb{N}_+$ by \\

(4.1)~~~ $ \Sigma f ( n ) ~=~ f ( 1 ) + f ( 2 ) + f ( 3 ) + \ldots + f ( n ),~~~
                                                               n \in \mathbb{N}_+ $ \\

Let \\

(4.2)~~~ $ \Phi $ \\

be a given set of functions $f : \mathbb{N}_+ \longrightarrow \mathbb{N}_+$. We denote by \\

(4.3)~~~ $ {\cal S}\Phi $ \\

its subset of all functions $f \in \Phi$, such that \\

(4.4)~~~ $ \Sigma f \in \Phi $ \\

and we are looking for {\it characterizations} of functions in ${\cal S}\Phi$. \\

We can further generalize the above scheme as follows. Given any family of integers \\

(4.5)~~~ $ \alpha = ( \alpha_1, \alpha_2, \alpha_3, \ldots ) \in \mathbb{Z}^{\mathbb{N}_+} $ \\

we extend the definition in (4.1) by \\

(4.6)~~~ $ \Sigma_\alpha f ( n ) ~=~ \alpha_1 f ( 1 ) + \alpha_2 f ( 2 ) +
                                       \alpha_3 f ( 3 ) + \ldots + \alpha_n f ( n ),~~~
                                                                    n \in \mathbb{N}_+ $ \\

Then, instead of (4.3), we consider the functions in \\

(4.7)~~~ $ {\cal S}_\alpha \Phi $ \\

which is the subset of all functions $f \in \Phi$, such that \\

(4.4)~~~ $ \Sigma_\alpha f \in \Phi $ \\

and now, we are looking for {\it characterizations} of functions in ${\cal S}_\alpha
\Phi$. \\ \\

\newpage

{\bf Appendix} \\

{\bf Lemma A} \\

Given $k \in \mathbb{N}_+$, then for $ \in \mathbb{N}_+$, we have \\

(A.1)~~~ $ 1^k + 2^k + 3^k + \ldots + n^k ~=~
                c_{k,\,0}\, n + c_{k,\, 1}\, n^2 +
                                      \ldots + c_{k,\, k}\, n^{k + 1} $ \\

for suitable $c_{k,\,0},~ c_{k,\, 1},~ \ldots ,~ c_{k,\, k} \in \mathbb{Q}$ which do not
depend on $n$. \\

{\bf Proof.} \\

We note that (A.1) is in fact an immediate consequence of what is known as Faulhaber's formula,
see http://mathworld.wolfram.com/ FaulhabersFormula.html

\end{document}